\documentclass[12pt,letterpaper,reqno,abstract=true]{amsart}

\usepackage{amsmath}
\usepackage{amsfonts}
\usepackage{amssymb}
\usepackage{amsthm}
\usepackage{amscd}

\usepackage{tikz}
\usepackage{tikz-cd}
\usetikzlibrary{cd}

\usepackage[colorlinks=true]{hyperref}
\hypersetup{urlcolor=red, citecolor=blue}

\usepackage{latexsym}
\usepackage{mathrsfs}
\usepackage{psfrag}
\usepackage[dvips]{epsfig}
\usepackage{epsfig}
\usepackage[all]{xy}         
\usepackage{wrapfig}

\usepackage[margin=2.5cm]{geometry}

\swapnumbers 

\theoremstyle{plain}                              

\newtheorem{thm}{Theorem}[section]

\newtheorem{defn}[thm]{Definition}

\newtheorem{cor}[thm]{Corollary}

\theoremstyle{definition}                         
\newtheorem{example}[thm]{Example}
\newtheorem{exercise}[thm]{Exercise}
\newtheorem*{remark}{Remark} 

\theoremstyle{remark}                             

\newcommand{\R}{\mathbb{R}}                     

\newcommand{\veps}{\varepsilon}

\providecommand{\abs}[1]{\left\lvert#1\right\rvert}        

\begin{document}

\title[Geometric control theory]{A short introduction to
  geometric control theory}

\author[S. N. Simi\'c]{Slobodan N. Simi\'c}

\address{Department of Mathematics and
  Statistics \\San Jos\'e State University \\ San Jose, CA 95192-0103}

\email{slobodan.simic@sjsu.edu}

\thanks{This work was partially supported by an SJSU
Research, Scholarship, and Creative Activity grant.}

\maketitle


\begin{abstract}
  The goal of this expository paper is to present the basics of
  geometric control theory suitable for advanced undergraduate or
  beginning graduate students with a solid background in advanced
  calculus and ordinary differential equations.
\end{abstract}






This paper consists of three parts:

\begin{enumerate}
\item Basic concepts;
\item Basic results;
\item Steering with piecewise constant inputs.
\end{enumerate}

The goal is to present only the necessary minimum to understand part
3, which describes a constructive procedure for steering affine
drift-free systems using piecewise constant inputs. All technical
details are omitted. The reader is referred to the literature at
the end for proofs and details.

\section{Basic concepts}
\label{sec:basic-concepts}

Control theory studies families of ordinary differential equations
parametrized by input, which is external to the ODEs and can be
controlled.

\begin{defn}
  A \textsf{control system} is an ODE of the form
  \begin{displaymath}
    \dot{x} = f(x,u),
  \end{displaymath}
  where $x \in M$ is the \emph{state} of the system, $M$ is the
  \emph{state space}, $u \in U(x)$ is the \emph{input} or
  \emph{control}, $U(x)$ is the (state-dependent) \emph{input set},
  and $f$ is a smooth map called the \emph{system map}. 
\end{defn}

\begin{figure}[h]
  \centering
  \includegraphics[scale=0.9]{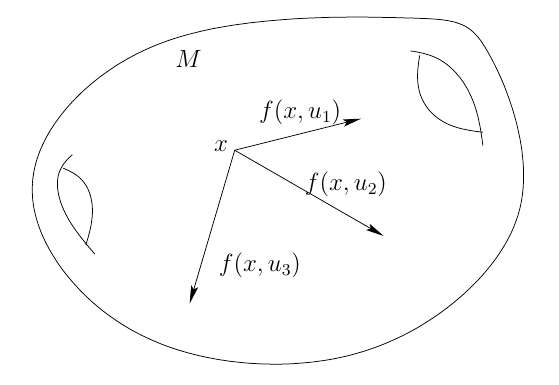}
  \caption{Control directions at a point.}
  \label{fig:1}
\end{figure}

The state space is usually a Euclidean space or a smooth manifold. The
set 
\begin{displaymath}
  \mathbb{U} = \bigcup_{x \in M} U(x)
\end{displaymath}
is called the \emph{control bundle}. With each control system we
associate the set of \emph{admissible control functions},
$\mathcal{U}$, consisting of functions $u : [0,T] \to \mathbb{U}$, for
some $T > 0$. These are usually square integrable functions, such as
piecewise continuous, piecewise smooth, or smooth ones (depending on
what the control system is modeling).

\begin{defn}
  A curve $x: [0,T] \to M$ is called a \textsf{control trajectory} if
  there exists an admissible control function
  $u : [0,T] \to \mathbb{U}$ such that $u(t) \in U(x(t))$ and
  \begin{displaymath}
    \dot{x}(t) = f(x(t),u(t)),
  \end{displaymath}
  for all $0 \leq t \leq T$.
\end{defn}

The basic idea is: if we want to move according to control system
$\dot{x} = f(x,u)$ starting from a point $x \in M$, the directions
that we have at our disposal are given by $f(x,u)$, for all
$u \in U(x)$. A control trajectory $x(t)$ defined by an admissible
control $u(t)$ picks one such direction for each time $t \in [0,T]$;
this choice is, of course, determined by $u(t)$.

\begin{figure}[h]
  \centering
  \includegraphics[scale=0.9]{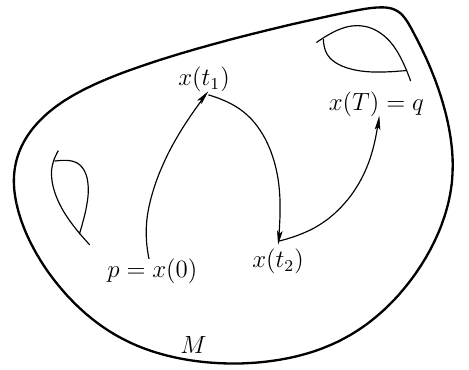}
  \caption{A control trajectory connecting $p$ and $q$.}
  \label{fig:2}
\end{figure}

\begin{example}[Linear control systems]
  If $M = \R^n$ and $f$ is linear in both $x$ and $u$, the system is
  called \emph{linear}:
  \begin{displaymath}
    \dot{x} = Ax + \sum_{i=1}^m u_i b_i = Ax + Bu,
  \end{displaymath}
  where $b_1,\ldots,b_m \in \R^n$ and $B$ is the $n \times m$ matrix
  whose columns are $b_1, \ldots, b_m$.
\end{example}

\begin{example}[Affine control systems]
  If $f$ is affine in $u$, the system is called \emph{affine}:
  \begin{displaymath}
    \dot{x} = g_0(x) + \sum_{i=1}^m u_i g_i(x),
  \end{displaymath}
  where $g_0,\ldots,g_m$ are vector fields on $M$; $g_0$ is called
  the \emph{drift}, and $g_1,\ldots, g_m$ are \emph{control vector
    fields}. If $g_0 = 0$, the system is called \emph{drift-free}.
\end{example}

The notions of reachability and controllability, fundamental to
control theory, are defined next.

\begin{defn}
  If $x : [0,T] \to M$ is a control trajectory with $x(0) = p$ and
  $x(T) = q$, we say that $q$ is \textsf{reachable} or
  \textsf{accessible} from $p$. The set of points reachable from $p$
  will be denoted by $\mathcal{R}(p)$.
\end{defn}

If the interior (relative to $M$) of $\mathcal{R}(p)$ is non-empty, we
say that the system is \textsf{locally accessible at} $p$. If it is
locally accessible at \emph{every} $p$, it is called \textsf{locally
  accessible}.

If $\mathcal{R}(p) = M$ for some (and therefore all) $p$, the system
is called \textsf{controllable}.

\begin{example}
  The system $\dot{x}_1 = u_1$, $\dot{x}_2 = u_2$, with $(u_1,u_2) \in
  \R^2$ (the inputs are \emph{unconstrained}) is trivially
  controllable.  \qed
\end{example}

\begin{example}
  Consider the control system $\dot{x} = \alpha x + u$, where
  $\alpha > 0$ and $u \in [-1,1]$ (the input is
  \emph{constrained}). We claim that it is \emph{uncontrollable}, for
  any choice of admissible control functions $\mathcal{U}$.

  To prove this, let $u \in \mathcal{U}$ be arbitrary and let $x(t)$
  be the corresponding control trajectory. The basic theory of linear
  ODEs yields
  \begin{align*}
    x(t) & = e^{\alpha t} x(0) + e^{\alpha t} \int_0^t e^{-\alpha s}
           u(s) \: ds \\
         & \geq x(0) - e^{\alpha t} \int_0^t e^{-\alpha s} \: ds \\
         & = x(0) + \frac{e^{\alpha t} - 1}{\alpha} \\
    & > x(0),
  \end{align*}
  for all $t > 0$. Thus any $q < p = x(0)$ is unreachable from $p$, so
  the system is uncontrollable. \qed
\end{example}

\begin{exercise}
  Consider the system
  \begin{align*}
    \dot{x}_1 & = x_2 \\
    \dot{x}_2 & = -kx_1-ux_2,
  \end{align*}
  where $k > 0$ is a constant and $u \in \R$ is the input. Show that
  the system is controllable on $\R^2 \setminus \{(0,0)\}$. More
  precisely, show that for every $p, q \in \R^2 \setminus \{(0,0)\}$,
  $q$ can be reached from $p$ using piecewise constant input with at
  most one switch between control vector fields. (This type of system
  occurs in mechanics.)
\end{exercise}

\subsection*{A geometric point of view}
Suppose we have an affine drift-free system
\begin{displaymath}
  \dot{x} = u_1 X_1(x) + \cdots + u_m X_m(x),
\end{displaymath}
where the inputs $u_1, \ldots, u_m \in \R$ are unconstrained. Then at
each point $x \in M$, the set of directions along which the system can
evolve is given by
\begin{displaymath}
  \Delta(x) = \text{span}\{ X_1(x), \ldots, X_m(x) \}.
\end{displaymath}
Note that $\Delta(x)$ is a subspace of $T_x M$, the tangent space to
$M$ at $x$. Recall that the tangent bundle of $M$ is defined by
\begin{displaymath}
  TM = \bigcup_{x \in M} T_x M.
\end{displaymath}
For instance, if $M$ is $\R^n$ or a Lie group of dimension $n$, then
$TM = M \times \R^n$. (This is not true in general; take for instance
$M = S^2$, the 2-sphere in $\R^3$.) Therefore, the evolution of the
control system is specified by a collection of planes $\Delta(x)$,
with $x \in M$. Such an object is called a \textsf{plane field} or a
\textsf{distribution}.

\begin{defn}
  A \textsf{distribution} on a smooth manifold $M$ is an assignment to each
  point $x \in M$ of a linear subspace $\Delta(x)$ of the tangent
  space $T_x M$.
\end{defn}

The distribution $\Delta$ associated with an affine drift-free control
system as above is called the \textsf{control distribution}.

Assume for a moment that $M = \R^n$, $m = 1$, and $X_1(x) \neq 0$, for
all $x$. Then $\Delta$ is 1-dimensional and through each point $x \in
\R^n$ there passes a curve $\mathcal{F}(x)$ (namely, the integral
curve of $X_1$) everywhere tangent to $\Delta$ (i.e., to $X_1$). This
collection of curves $\mathcal{F}(x)$ is called a \textsf{foliation}
which integrates $\Delta$. It is not hard to see that
\begin{displaymath}
  \mathcal{F}(x) = \mathcal{R}(x),
\end{displaymath}
that is, $\mathcal{F}(x)$ is exactly the set of points reachable from
$x$. Since $\mathcal{F}(x) \neq \R^n$, the system is not
controllable. Thus at least in this simplified case, we have
\begin{displaymath}
  \Delta \ \text{is integrable} \quad \Rightarrow \quad \text{system is uncontrollable}.
\end{displaymath}
So if the system is controllable, its control distribution should
satisfy a property that is, intuitively, opposite to integrability. We
will soon see that $\Delta$, in fact, has to be \emph{bracket
  generating}.

\section{Basic results}
\label{sec:basic-results}

For a smooth vector field $X$ on $M$, denote by $\{ X^t \}$ its
(local) \emph{flow}. That is, for each $p \in M$, $t \mapsto X^t(p)$
is the integral curve of $X$ staring at $p$. Recall that $X^{s+t} =
X^s \circ X^t$, wherever both sides are defined. If $M$ is compact or
if $X$ is bounded, then $X^t$ is defined for all $t$, i.e., the flow
is \emph{complete}. In that case $X^t : M \to M$ is a diffeomorphism
(i.e., smooth together with its inverse).

For a diffeomorphism $\phi : M \to M$, denote by $\phi_\ast$ the
associated \textsf{push-forward map} defined on vector fields by
\begin{displaymath}
  \phi_\ast(X)(p) = (T_{\phi^{-1}(p)} \phi)(X(\phi^{-1}(p))),
\end{displaymath}
where $X$ is a vector field on $M$ and $T_q\phi$ denotes the derivative
(or tangent map) of $\phi$ at $q$.

\begin{defn}
  The \textsf{Lie bracket} of smooth vector fields $X, Y$ is defined
  by
  \begin{displaymath}
    [X,Y] = \left. \frac{d}{dt} \right|_{t=0} (X^{-t})_\ast(Y).
  \end{displaymath}
\end{defn}

That is, $[X,Y]$ is the derivative of $Y$ along the integral curves of
$X$. It turns out that $[X,Y]$ can also be expressed in the following
way:
\begin{displaymath}
  [X,Y](p) = DY(p)X(p) - DX(p) Y(p).
\end{displaymath}
For diffeomorphism $\phi, \psi : M \to M$ define their
\textsf{bracket} by
\begin{displaymath}
  [\phi,\psi] = \psi^{-1} \circ \phi^{-1} \circ \psi \circ \phi.
\end{displaymath}

\begin{thm}[The fundamental fact about Lie brackets] \label{thm:fundamental}
  If $[X,Y] = Z$, then
  \begin{displaymath}
    [X^t,Y^t](p) = Z^{t^2}(p) + o(t^2),
  \end{displaymath}
  as $t \to 0$. 
\end{thm}

The Landau ``little $o$'' notation $f(t) = o(t^2)$ means that $f(t)/t^2 \to
0$, as $t \to 0$.

\begin{figure}[h]
  \centering
  \includegraphics[scale=0.9]{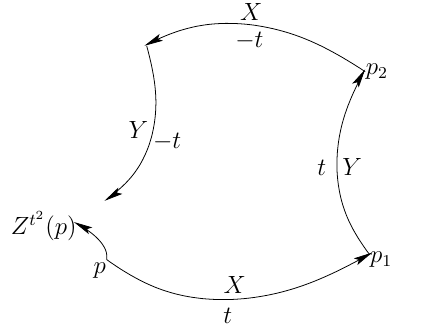}
  \caption{The fundamental fact about the Lie bracket.}
  \label{fig:3}
\end{figure}

Therefore, to advance in the direction of the Lie bracket $Z = [X,Y]$
while moving only along integral curves of $X$ and $Y$, we need to
move from $p$ to $p_1 = X^t(p)$ (along $X$), then to $p_2 = Y^t(p_1)$
(along $Y$), etc. (Figure~\ref{fig:3}). From the point of view of
control theory, this is crucial information if
$Z \not\in \text{span} \{X,Y \}$.

\begin{remark}
  Observe that the “price to pay” to move $t^2$ units of time in the
  $Z$-direction using only $X$ and $Y$ is $4t$ units of time. When $t$
  is very small, $4t$ is much larger than $t^2$.
\end{remark}

\begin{defn}
  A distribution $\Delta$ is called \textsf{involutive} if for every
  two vectors fields $X, Y$ in $\Delta$ (in the sense that
  $X(p), Y(p) \in \Delta(p)$), $[X,Y]$ is in $\Delta$.  
\end{defn}

That is, $\Delta$ is involutive if it is closed under the Lie bracket.

A fundamental result in the theory of smooth manifolds is the
following:

\begin{thm}[Frobenius's theorem]
  If a smooth distribution $\Delta$ on $M$ of constant dimension $k$
  is involutive, then it is \textsf{completely integrable}. That is,
  through every point $p \in M$ there passes a $k$-dimensional
  (immersed) submanifold $\mathcal{F}(p)$ (called an \textsf{integral
    manifold} of $\Delta$) which is everywhere tangent to $\Delta$;
  i.e., for every $q \in \mathcal{F}(p)$,
  $T_q \mathcal{F}(p) = \Delta(q)$. The collection $\mathcal{F}(p)$
  ($p \in M$) forms the \textsf{integral foliation} of $\Delta$.

  Equivalently, for every $p \in M$ there is a neighborhood $U$ of $p$
  and a diffeomorphism $\phi : U \to \phi(U)$, with
  $\phi(q) = (x_1(q), \ldots, x_n(q))$ such that the slices
  $x_{k+1} = \text{\rm constant}, \ldots, x_n = \text{\rm constant}$
  is an integral manifold of $\Delta$.
\end{thm}

If $\Delta$ happens to be the control distribution of an affine
drift-free system, then for every $p \in M$, we have $\mathcal{R}(p) =
\mathcal{F}(p)$. That is, the set of points reachable from $p$ is
precisely the integral manifold $\mathcal{F}(p)$.

Thus if we want a control system to be controllable, its control
distribution needs to have a property that is diametrically opposite
to involutivity. That property is defined as follows.

\begin{defn}
  A distribution $\Delta = \text{span} \{ X_1, \ldots, X_k \}$ on $M$
  is \textsf{bracket generating} if the iterated Lie brackets
  \begin{displaymath}
    X_i, [X_i,X_j], [X_i,[X_j,X_k]], \ldots
  \end{displaymath}
  with $1 \leq i, j, k, \ldots \leq m$ span the tangent space of $M$
  at every point. 
\end{defn}

In other words, one can obtain all directions in the tangent space to
$M$ by taking iterated Lie brackets of vector fields in $\Delta$. It
turns out that the property of being bracket generating does not
depend on the choice of the frame $X_1, \ldots, X_m$.

The fundamental result in geometric contro theory is the following theorem.

\begin{thm}[Chow-Rashevskii]
  If $\Delta =\text{\rm span}\{X_1,\ldots, X_m \}$ is a bracket generating
  distribution on $M$, then every two points in $M$ can be connected
  by a path which is almost everywhere tangent to $\Delta$. The path
  can be chosen to be piecewise smooth, consisting of arcs of integral
  curves of $X_1,\ldots, X_m$.
\end{thm}

In the language of control theory, we have:

\begin{cor}
  If the control distribution of an affine drift-free control system
  is bracket generating, then the system is controllable.
\end{cor}

\begin{example}[The Heisenberg group]
  On $\R^3$ define a distribution $\Delta = \text{span}\{ X_1, X_2
  \}$, where
  \begin{displaymath}
    X_1(x,y,z) =
    \begin{bmatrix}
      1 \\ 0 \\ 0
    \end{bmatrix}, \qquad
    X_2(x,y,z) =
    \begin{bmatrix}
      0 \\ 1 \\ x
    \end{bmatrix}.
  \end{displaymath}
  Observe that
  \begin{displaymath}
    X_3 := [X_1,X_2] =
    \begin{bmatrix}
      0 \\ 0 \\ 1
    \end{bmatrix}.
  \end{displaymath}
  Since $X_1, X_2, X_3$ span the tangent space to $\R^3$ at every
  point, $\Delta$ is a 2-dimensional bracket generating
  distribution. The triple $(\R^3, \Delta, \langle \cdot, \cdot
  \rangle)$, where $\langle \cdot, \cdot \rangle$ is the restriction
  of the usual dot product to $\Delta$, is called the
  \textsf{Heisenberg group}, a basic example of a \emph{subriemannian geometry}.

  Note that at the origin, $X_1(0,0,0) = (1,0,0)^T$ and $X_2(0,0,0) =
  (0,1,0)^T$. It is natural to ask: how to we reach the point
  $(0,0,z)$ (with, e.g., $z > 0$) from the origin by a control
  trajectory of the control system $\dot{x} = u_1 X_1(x) + u_2
  X_2(x)$?

  Theorem~\ref{thm:fundamental} and $[X_1,X_3] = [X_2,X_3] = 0$ imply
  that
  \begin{displaymath}
    X_3^{t^2} = [X_1^t,X_2^t],
  \end{displaymath}
  with exact rather than asymptotic equality. It follows that for $z >
  0$, we have
  \begin{displaymath}
    (0,0,z) = [X_1^{\sqrt{z}},X_2^{\sqrt{z}}](0,0,0).
  \end{displaymath}
  Therefore,a piecewise constant input function which steers the
  control system $\dot{x} = u_1 X_1(x) + u_2 X_2(x)$ from $(0,0,0)$ to
  $(0,0,z)$ (with $z > 0$) is the following:
  \begin{displaymath}
    u(t) = (u_1(t),u_2(t) = 
    \begin{cases}
      (1,0), & \text{for} \quad 0 \leq t < \sqrt{z} \\
      (0,1), & \text{for} \quad \sqrt{z} \leq t < 2 \sqrt{z} \\
      (-1,0), & \text{for} \quad 2 \sqrt{z} \leq t < 3 \sqrt{z} \\
      (0,-1), & \text{for} \quad 3 \sqrt{z} \leq t < 4 \sqrt{z}.
    \end{cases}
  \end{displaymath}
  If $z < 0$, we can simply reverse the order of $X_1$ and $X_2$, and
  use $\sqrt{-z}$ instead of $\sqrt{z}$. \qed  
\end{example}

\subsection*{A brief look into sub-Riemannian geometry}
In the Heisenberg group, for any two points $p, q \in \R^3$ there
exists a \emph{horizontal path} $c$ (i.e., a path almost everywhere
tangent to $\Delta$) which starts at $p$ and terminates at $q$. It is
natural to ask, what is the length of the shortest such path? Let
\begin{displaymath}
  d_\Delta(p,q) = \inf \{ \ell(c) : c \ \text{is a horizontal path
    and} \ c(0) = p, c(1) = q \},
\end{displaymath}
where $\ell(c)$ os the arc-length of $c$. This well-defined number is
called the \textsf{sub-Riemannian distance} between $p$ and $q$. It
turns out that for every two points $p$ and $q$ there exists a
horizontal path $c_0$, called a \textsf{sub-Riemannian geodesic}, which
realizes this distance: $d_\Delta(p,q) = \ell(c_0)$. Sub-Riemannian
geometry is closely related relative of \textsf{optimal control
  theory}, which is beyond the scope of this paper.

\section{Steering affine drift-free systems using piecewise constant
  inputs}
\label{sec:steer-affine-drift}

The goal of this section is to outline an answer to the question:

\begin{quote}
  \textit{How does one steer a general affine drift-free system?}
\end{quote}

We begin by recalling that in the Heisenberg group, all Lie brackets
of the generating vector fields $X_1, X_2$ of order $> 2$ vanish. We
say that the Lie algebra generated by $X_1, X_2$ is \textsf{nilpotent
  of order 2}.

\begin{defn}
  A Lie algebra $\mathcal{L}$ is called \textsf{nilpotent of order} $k$
  if all Lie brackets of elements of $\mathcal{L}$ of order $> k$ vanish.
\end{defn}

Now consider an affine control system on an $n$-dimensional manifold
$M$:
\begin{equation}    \label{eq:2}
  \dot{x} = u_1 X_1(x) + \cdots + u_m X_m(x).
\end{equation}
Let $\Delta$ be its control distribution. Assume that

\begin{enumerate}
\item $\Delta$ is bracket generating, and
\item the Lie algebra $\mathcal{L}$ generated by $X_1, \ldots, X_m$ is
  nilpotent of order $k$.
\end{enumerate}

To make the notation more intuitive, we will denote the flow of a
vector field $X$ by $e^{tX}$ , and will pretend that
\begin{displaymath}
  e^{tX} = I + tX + \frac{t^2}{2!} X^2 + \frac{t^3}{3!} X^3 + \cdots,
\end{displaymath}
which of course makes sense and is correct for \emph{linear} vector
fields. In general, this expression makes no sense without further
clarification, but we can still treat it as a formal series in some
“free Lie algebra” (namely $\mathcal{L}$, if $X = X_i$, for some $i$).

As a vector space, the Lie algebra $\mathcal{L}$ admits a basis
$B_1,\ldots, B_s$ called the \textsf{Philip Hall basis}. This is a
canonically chosen basis of $\mathcal{L}$, which takes into account
the Jacobi identity. (We won’t go into details of how to compute this
basis.) The vector fields $B_i$ are suitably chosen Lie brackets of
$X_1, \ldots, X_m$, with $B_i = X_i$, for $1 \leq i \leq m$. (Thus
possibly $s > n$.)

The so called \textsf{Chen-Fliess formula} asserts us that every
“flow” (i.e., control trajectory) of \eqref{eq:2} is of the form
\begin{equation}  \label{eq:CF}
  S(t) = e^{h_s(t) B_s} \cdots e^{h_1(t) B_1},
\end{equation}
for some real-valued functions $h_1, \ldots, h_s$ called the
\textsf{Philip Hall coordinates}. Furthermore, we have $S(0) = I$ (the
identity) and
\begin{equation}   \label{eq:4}
  \dot{S}(t) = S(t) \{ v_1(t) B_1 + \cdots + v_s(t) B_s \},
\end{equation}
where $v_1(t), \ldots, v_s(t)$ are called \textsf{fictitious
  inputs}. (Why fictitious? Only $v_1,\ldots, v_m$ are “real”,
corresponding to the vector fields $X_1,\ldots, X_m$.)

Let $p, q \in M$ be arbitrary. The algorithm due to Sussmann and
Lafferriere for steering the control system \eqref{eq:2} from $p$ to
$q$ is defined as follows:

\begin{enumerate}
\item[{\sc Step 1}] Find fictitious inputs steering the \emph{extended
    system}
  \begin{displaymath}
    \dot{x} = v_1 B_1(x) + \cdots + v_s B_s(x)
  \end{displaymath}
  from $p$ to $q$.

\item[{\sc Step 2}] Find “real” inputs $u_1,\ldots, u_m$ which
  generate the same evolution as $v_1,\ldots , v_s$.
  
\end{enumerate}

Details follow.

\medskip

\noindent \textsc{Step 1:} First note that by bracket generating
property of $\Delta$, we have $s \geq n = \dim M$. Thus the first step
is easy: simply take any curve $\gamma$ connecting $p$ and $q$, and
for each $t$, expres $\dot{\gamma}(t)$ as a linear combination of
$B_1, \ldots, B_s$. The coefficients are the desired fictitious
inputs:
\begin{displaymath}
  \dot{\gamma}(t) = \sum_{i=1}^s v_i(s) B_i(\gamma(t)).
\end{displaymath}
\noindent \textsc{Step 2:} Next, differentiate \eqref{eq:CF} with
respect to $t$ and use the chain rule:
\begin{align}
  \dot{S}(t) & = \sum_{i=1}^s e^{h_s(t) B_s} \cdots e^{h_{i+1}
               B_{i+1}} \cdot \dot{h}_i(t) B_i \cdot e^{h_{i-1}(t)
               B_{i-1}} \cdots e^{h_1(t) B_1} \nonumber \\
             & = \sum_{i=1}^s S(t) S_i(t) \{ \dot{h_i}(t) B_i \}
               S_i^{-1}(t) \nonumber \\
             & = \sum_{i=1}^s S(t) \text{Ad}_{S_i(t)}(\dot{h_i}(t)
               B_i),  \label{eq:*}
\end{align}
where
\begin{displaymath}
  S_i(t) = e^{-h_1(t) B_1} \cdots e^{-h_{i-1} B_{i-1}}
\end{displaymath}
and
\begin{displaymath}
  \text{Ad}_T(X) = T X T^{-1}.
\end{displaymath}
Since $\mathcal{L}$ is nilpotent, it follows that
\begin{displaymath}
  \text{Ad}_{S_i(t)}(B_i) = \sum_{i=1}^s p_{ij}(h(t)) B_j,
\end{displaymath}
for some \emph{polynomials} $p_{ij}(h_1,\ldots, h_s)$, where $h(t) =
(h_1(t),\ldots,h_s(t))$. Substituting into \eqref{eq:*}, we obtain
\begin{displaymath}
  \dot{S}(t) = S(t) \sum_{j=1}^s \left\{  \sum_{i=1}^s p_{ij}(h(t))
    \dot{h}_i(t)  \right\} B_j.
\end{displaymath}
A comparison with \eqref{eq:4} yields the following system of equation
for the fictitious inputs $v = (v_1,\ldots,v_s)^T$:
\begin{displaymath}
  v(t) = P(h(t)) \dot{h}(t),
\end{displaymath}
where $P(h) = [p_{ij}(h)]_{1 \leq i, j \leq s}$. It can be shown that
$P(h)$ is invertible; denote its inverse by $Q(h)$. Then
\begin{displaymath}
  \dot{h} = Q(v) v.
\end{displaymath}
This is called the \textsf{Chen-Fliess-Sussmann equation}. Given $v$,
obtained in {\sc Step 1}, we can use it to solve for $h$.

What remains to be done is to find ``real'' piecewise constant inputs
$u_1, \ldots , u_m$ which generate the same motion. The basic idea is
to use the Theorem~\ref{thm:fundamental} and the
\emph{Baker-Campbell-Hausdorff} formula to express the term $o(t^2)$ in
terms of the higher order Lie brackets. Instead of showing how to
do this in general, here is an example. Here's an example of how to do
this.

\begin{example}
  Consider the case $n = 4, m = 2$, and $k = 3$. This corresponds to a
  control system $\dot{x} = u_1 X_1(x) + u_2 X_2(x)$ in $\R^4$ such
  that all Lie brackets of order $>3$ vanish.

  The Philip Hall basis of $\mathcal{L}$, the Lie algebra generated by
  $X_1, X_2$, is
  \begin{align*}
    B_1 & = X_1 \\
    B_2 & = X_2 \\
    B_3 & = [X_1,X_2] \\
    B_4 & = [X_1,[X_1,X_3]] =
          [B_1,B_2] \\
    B_5 & = [X_2,[X_1,X_2]] = [B_2,B_3].
  \end{align*}
  If we think of $B_5$ as a vector field in $\R^4$ (which we will need
  to do eventually), then $B_5$ can be expressed as a linear
  combination of $B_1, \ldots, B_4$, so we can take $v_5(t) =
  0$. However, note that we cannot disregard $B_5$, because we need a
  complete basis of $\mathcal{L}$ in order for the algorithm to work.

  Setting $s = 5$, differentiating \eqref{eq:CF}, and expressing
  everything in terms of $B_1, \ldots, B_5$, we obtain:
  \begin{align*}
    \dot{h}_1(t) & \quad \text{multiplies} \ B_1  \\
    \dot{h}_2(t) & \quad \text{multiplies} \ B_2 - h_1 B_3 +\frac{h_1^2}{2}
                   B_4 \\
    \dot{h}_3(t) & \quad \text{multiplies} \ B_3 - h_2 B_5 - h_1 B_4 \\
    \dot{h}_4(t) & \quad \text{multiplies} \ B_4 \\
    \dot{h}_5(t) & \quad \text{multiplies} \ B_5.
  \end{align*}
  The Chen-Fliess-Sussmann equation is therefore
  \begin{align*}
    \dot{h}_1 & = v_1 \\
    \dot{h}_2 & = v_2 \\
    \dot{h}_3 & = h_1 v_2 + v_3 \\
    \dot{h}_4 & = \frac{h_1^2}{2} v_2 + h_1 v_3 + v_4 \\
    \dot{h}_5 & = h_2 v_3 + h_1 h_2 v_2.
  \end{align*}
  Let us now assume that we would like to generate
  \begin{displaymath}
    S(T) = e^{\veps B_5} e^{\delta B_4} e^{\gamma B_3} e^{\beta B_2} e^{\alpha B_1},
  \end{displaymath}
  for some given numbers $\alpha, \beta, \gamma, \delta, \veps$.

  Remember that we only have \emph{two inputs} available to achieve
  this. Denote by $w_i$ the input that gives rise to $e^{X_i}$. It is
  easy to see that $w_1 = (1,0)$ and $w_2 = (0,1)$. It follows that
  $\alpha w_1, \beta w_2$ generate $e^{\alpha B_1}, e^{\beta B_2}$,
  respectively. Denote the \emph{concatenation} of paths by the symbol
  $\sharp$. Then
  \begin{displaymath}
   \alpha w_1 \sharp \beta w_2
  \end{displaymath}
  generates $e^{\beta B_2} e^{\alpha B_1}$.

  By the Baker-Campbell-Hausdorff formula and nilpotency, it follows
  that
  \begin{displaymath}
    \sqrt{\gamma} w_1 \sharp \sqrt{\gamma} w_2 \sharp (-\sqrt{\gamma}
    w_1) \sharp (-\sqrt{\gamma} w_2)
  \end{displaymath}
  generates
  \begin{displaymath}
    e^{\hat{\gamma} B_5} e^{\hat{\gamma} B_4} e^{\gamma B_3},
  \end{displaymath}
  where $\hat{\gamma} = \frac{1}{2} \gamma^{3/2}$. It remains to
  generate $e^{\veps B_5} e^{\delta B_4} e^{-\hat{\gamma} B_4}
  e^{\hat{\gamma} B_5}$.

  Let $\rho = (\delta - \hat{\gamma})^{1/3}$ and $\sigma = (\veps +
  \hat{\gamma})^{1/3}$. A long and tedious calculation (which we omit)
  shows that
  \begin{eqnarray}
    \rho w_1 \sharp \rho w_1 \sharp \rho w_2 \sharp (-\rho w_1) \sharp
    (-\rho w_2) \sharp (-\rho w_1) \sharp \rho w_2 \sharp \rho w_1
    \sharp(-\rho w_2) \sharp (-\rho w_1) \\
    \sharp \sigma w_2 \sharp \sigma w_1 \sharp \sigma w_2 \sharp (-\sigma w_1) \sharp
    (-\sigma w_2) \sharp (-\sigma w_2) \sharp \sigma w_2 \sharp \sigma w_1
    \sharp(-\sigma w_2) \sharp (-\sigma w_1).
  \end{eqnarray}
  Concatenating the above four pieces together, we obtain a control
  consisting of \textbf{26 pieces}. \qed
\end{example}

What if the Lie algebra generated by the control vector fields is not
nilpotent? Then there are at least two possibilities:

\begin{enumerate}
\item[(a)] We can steer approximately using the above algorithm, or
\item[(b)] We can try to \emph{nilpotentize} the system, i.e.,
  reparametrize it so that the new system (more precisely, the
  corresponding Lie algebra) becomes nilpotent.
\end{enumerate}

Here is an example how the latter.

\begin{example}[Model of a unicycle]
  Consider the following model of a unicycle:
  \begin{displaymath}
    \dot{p} = u_1 X_1(p) + u_2 X_2(p),
  \end{displaymath}
  with $p = (x,y,z)$, where $(x,y) \in \R^2$ is the \emph{position} of
  the unicycle and $z \in [0,2\pi]$ is the \emph{angle} between the
  wheel and the positive $x$-axis. The control vector fields are
  \begin{displaymath}
    X_1(p) =
    \begin{bmatrix}
      \cos z \\ \sin z \\ 0
    \end{bmatrix}
    \quad \text{and} \quad
    X_2(p) =
    \begin{bmatrix}
      0 \\ 0 \\ 1
    \end{bmatrix}.
  \end{displaymath}
  It is not hard to check that the Lie algebra generated by $X_1, X_2$
  is not nilpotent. However, if $\abs{z} < \pi/2$, we can
  reparametrize the system using the following \emph{feedback
    transformation}:
  \begin{displaymath}
    u_1 = \frac{1}{\cos z} v_1, \qquad u_2 = (\cos^2 z) v_2.
  \end{displaymath}
  We obtain a new system,
  \begin{displaymath}
    \dot{p} = v_1 Y_1(p) + v_2 Y_2(p),
  \end{displaymath}
  where
  \begin{displaymath}
    Y_1(p) =
    \begin{bmatrix}
      1 \\ \tan z \\ 0
    \end{bmatrix}
    \quad \text{and} \quad
    Y_2 =
    \begin{bmatrix}
      0 \\ 0 \\ \cos^2 z
    \end{bmatrix}.
  \end{displaymath}
  Since
  \begin{displaymath}
    [Y_1,Y_2] = Y_3 =
    \begin{bmatrix}
      0 \\ -1 \\ 0
    \end{bmatrix},
    \quad [Y_1,Y_3] = [Y_2,Y_3] = 0,
  \end{displaymath}
  the new system is indeed nilpotent of order 2. Furthermore, in the
  region defined by $\abs{z} < \pi/2$, $Y_1, Y_2$, and $Y_3$ span the
  tangent space. (Note that the feedback transformation did not change
  the bracket generating property of the control distribution.) \qed
\end{example}

\section{Further reading}
\label{sec:further-reading}

For a gentle and well-written introduction to smooth manifolds, I
recommend Boothby's \cite{boothby+2003}. The basics of nonlinear
control theory can be found in \cite{sastry+99} and
\cite{nijmeijer+90}. The former has more details on steering and an
extensive bibliography. For an introduction to geometric control, see
\cite{jurdj97}. An excellent introduction to sub-Riemannian geometry
is \cite{mont02}. For details on steering (that is, motion planning),
see \cite{laff+sussmann+93}.


\bibliographystyle{amsalpha} 





\end{document}